\documentclass[a4paper]{amsart}
\usepackage{latexsym,bm}
\usepackage{amsmath,amsfonts,amssymb,mathrsfs}
\usepackage{amscd}
\usepackage{cite}

\newcommand{\BS}{\mathfrak S}

\newcommand{\Z}{\mathbb Z}
\newcommand{\HH}{\mathscr{H}}

\newcommand{\Q}{\mathbb Q}
\newcommand{\lam}{\lambda}
\newcommand{\eps}{\varepsilon}

\newcommand{\ksl}{\widehat{\mathfrak{sl}}}
\newcommand\Lam{\Lambda}
\newcommand\umu{{\boldsymbol\mu}}
\newcommand\bmu{{\boldsymbol\mu}}

\newcommand{\bQ}{\mathbf{Q}}
\newcommand\blam{{\boldsymbol\lambda}}
\newcommand\ulam{{\boldsymbol\lambda}}
\newcommand\bnu{{\boldsymbol\nu}}
\newcommand\uempty{{\boldsymbol\emptyset}}

\DeclareMathOperator{\res}{res} \DeclareMathOperator{\wt}{wt}
\DeclareMathOperator{\rad}{rad} 
\renewcommand\t{\mathfrak{t}}
\DeclareMathOperator{\std}{Std}
\let\gedom=\trianglerighteq

\newcounter{main}
\theoremstyle{plain}

\swapnumbers
\newtheorem{prop}[equation]{Proposition}
\newtheorem{thm}[equation]{Theorem}
\newtheorem{cor}[equation]{Corollary}
\newtheorem{lem}[equation]{Lemma}
\newtheorem{example}[equation]{Example}
\newtheorem{dfn}[equation]{Definition}
\newtheorem{conj0}[equation]{Generalised DJM Conjecture}

\numberwithin{equation}{section}

{\catcode`\|=\active
  \gdef\set#1{\mathinner{\lbrace\,{\mathcode`\|"8000%
                                  \let|\midvert #1}\,\rbrace}}
}
\def\midvert{\egroup\mid\bgroup}

\title
  {On a generalisation of the Dipper--James--Murphy Conjecture}
\subjclass[2000]{20C08, 20C30} \keywords{Crystal basis, Fock spaces, Kleshchev multipartitions,
ladder multipartitions, ladder nodes, Lakshimibai--Seshadri paths}
\author{Jun Hu}
\address{School of Mathematics and Statistics\\
\newline
\null\,\quad University of Sydney\\
\newline\null\,\quad NSW 2006, Australia
\smallskip
\newline\null\,\qquad\qquad\qquad\&
\smallskip
\newline\null\,\quad
Department of Mathematics\\
\newline\null\,\quad Beijing Institute of Technology\\
\newline\null\,\quad Beijing 100081, P.R. China\\
\medskip} \email{junhu303@yahoo.com.cn}

\begin{document}
\bibliographystyle{andrew}

\begin{abstract}
Let $r, n$ be positive integers. Let $e$ be $0$ or an integer bigger
than $1$. Let $v_1,\cdots,v_r\in\Z/e\Z$ and $\mathcal{K}_r(n)$ the
set of Kleshchev $r$-partitions of $n$ with respect to $(e;\bQ)$,
where $\bQ:=(v_1,\cdots,v_r)$. The Dipper--James--Murphy conjecture
asserts that $\mathcal{K}_r(n)$ is the same as the set of
$(\bQ,e)$-restricted bipartitions of $n$ if $r=2$. In this paper we
consider an extension of this conjecture to the case where $r>2$. We
prove that any multi-core $\ulam=(\lam^{(1)},\cdots,\lam^{(r)})$ in
$\mathcal{K}_r(n)$ is a $(\bQ,e)$-restricted $r$-partition. As a
consequence, we show that in the case $e=0$, $\mathcal{K}_r(n)$
coincides with the set of $(\bQ,e)$-restricted $r$-partitions of $n$
and also coincides with the set of ladder $r$-partitions of $n$.
\end{abstract}

\maketitle


\section{Introduction}

A composition $\alpha=(\alpha_1,\alpha_2,\cdots)$ is a finite
sequence of non-negative integers; we denote by $|\alpha|$ the sum
of this sequence and call $\alpha$ a composition of $|\alpha|$. A
partition is a composition whose parts are non-increasing. Let $r,n$
be positive integers. A multipartition, or $r$-partition, of $n$ is
an ordered sequence $\blam=(\lambda^{(1)},\dots,\lambda^{(r)})$ of
partitions such that $|\lambda^{(1)}|+\dots+|\lambda^{(r)}|=n$. The
partitions $\lambda^{(1)},\dots,\lambda^{(r)}$ are called the
components of $\blam$.  If $r=2$, a multipartition is also called a
bipartition.

Let $e$ be $0$ or an integer bigger than $1$. Let
$v_1,\cdots,v_r\in\Z/e\Z$. If $e>1$, then a partition
$\lam=(\lam_1,\lam_2,\cdots)$ of $k$ is said to be
$e$-restricted if $\lam_i-\lam_{i+1}<e$ for any $i\geq 1$. We make
the convention that every partition is $e$-restricted if $e=0$. The
notion of $e$-restricted partitions plays an important role in the
modular representation theory of the symmetric groups $\BS_n$ as
well as its associated Iwahori--Hecke algebra $\HH_q(\BS_n)$. For
example, if $e>1$ and the parameter $q$ is a primitive $e$-th root
of unity, then it is well-known that simple modules of
$\HH_q(\BS_n)$ are in one-to-one correspondence with the set of
$e$-restricted partitions of $n$. The same is true if $e=0$ and $q$
is not a root of unity. Another important application (cf.
\cite{MM}) is that the set of $e$-restricted partitions provides a
combinatorial realization of the crystal graph of the integrable
highest weight module of level one over the affine Lie algebra
$\widehat{\mathfrak{sl}}_e$ if $e>1$; or over the affine Lie algebra
${\mathfrak{gl}}_{\infty}$ if $e=0$.

In \cite{AM}, Ariki and Mathas introduced a notion of Kleshchev
multipartitions which provides a combinatorial realization of the
crystal graph of integrable highest weight module of level $r$ over
the affine Lie algebra $\widehat{\mathfrak{sl}}_e$ if $e>1$; or over
the affine Lie algebra ${\mathfrak{gl}}_{\infty}$ if $e=0$. {\it A
priori}, the notion of Kleshchev multipartition is defined with
respect to the given $(r+1)$-tuple $(e;v_1,\cdots,v_r)$ and is
recursively defined. It is desirable to look for a non-recursive
definition. In the case $r=1$, it coincides with the notion of
$e$-restricted partitions. In general, by a result of  Ariki
\cite{A3}, the notion of Kleshchev multipartitions fits nicely with
the Dipper--James--Mathas Specht module theory of the cyclotomic
Hecke algebra $\HH_{r,n}(q;q^{v_1},\cdots,q^{v_r})$ and gives
natural labelling of the simple modules of
$\HH_{r,n}(q;q^{v_1},\cdots,q^{v_r})$, where the parameter $q$ is a
primitive $e$-th root of unity if $e>1$; or not a root of unity if
$e=0$. Thus, the notion of Kleshchev multipartitions can be regarded as a natural
generalization of the notion of $e$-restricted partitions.

In 1995, when $r=2$, Dipper, James and Murphy (see \cite{DJM2})
proposed a notion of $(Q,e)$-restricted bipartitions (which is
non-recursively defined), where $Q=-q^{m}$, $q=\sqrt[e]{1}$ and
$e>1, m\in\Z/e\Z$, and they conjectured that a Kleshchev bipartition
of $n$ with respect to $(e;m,1)$ is the same as a $(Q,e)$-restricted
bipartition of $n$. This conjecture was proved only recently by
Ariki--Jacon \cite{AJ}, using the result of another recent work of
Ariki--Kreiman--Tsuchioka \cite{AKT}. The paper \cite{AKT} contains
a new non-recursive description of Kleshchev bipartitions. In
general, in the case $r>2$, the question of finding a non-recursive
characterization of Kleshchev $r$-partitions remains open.\smallskip

The starting point of this paper is to explore this open question.
We give a natural extension of the Dipper--James--Murphy notion of
$(Q,e)$-restricted bipartitions to the case where $r>2$, i.e.,
$(\bQ,e)$-restricted multipartitions, where $\bQ:=(v_1,\cdots,v_r)$.
We also introduce a notion of ladder $r$-partitions. It turns
out that any $(\bQ,e)$-restricted multipartition of $n$ is a
Kleshchev multipartition in $\mathcal{K}_r(n)$. Our main result
asserts that any multi-core $\ulam=(\lam^{(1)},\cdots,\lam^{(r)})$
in $\mathcal{K}_r(n)$ is a $(\bQ,e)$-restricted multipartition.

As a consequence, we show
that if $e=0$ (in that case every multipartition is a multi-core), then $\mathcal{K}_r(n)$ coincides with the set
of $(\bQ,e)$-restricted multipartitions of $n$, which gives a
non-recursive description of Kleshchev $r$-partition in this
case; and also coincides with the set of ladder $r$-partitions of $n$, which
gives a new recursive description of Kleshchev $r$-partition in that
case. The main result is a generalization of the theorem of Ariki and Jacon \cite{AJ}, i.e., we prove a generalization of
the Dipper--James--Murphy conjecture to the case where $e=0$ and $r>2$. Conjecturally, everything should be still true in the case where $e>1$.
\smallskip

The paper is organized as follows. In Section 2, we recall the
notions of Kleshchev multipartitions and $(\bQ,e)$-restricted
multipartitions. In particular, we show that any
$(\bQ,e)$-restricted $r$-partition of $n$ is a Kleshchev
$r$-partition with respect to $(e;\bQ)$. We also recall a result of
Littelmann and a related result of Kashiwara, and give some
consequence of these two results. In Section 3, after introducing
the notion of ladder nodes, ladder sequences, ladder multipartitions
as well as strong ladder multipartitions, we give the proof of our
main result Propsition \ref{corecase}. As a consequence we prove the
generalized Dipper--James--Murphy conjecture when $e=0$, where we
also show that the notion of ladder $r$-partition coincides
with the notion of strong ladder $r$-partition in that case.

\bigskip
\section{Preliminaries}

Let $r, n$ be positive integers. Let $\mathcal{P}_{r}(n)$ be the set
of $r$-partitions of $n$. If $\ulam\in\mathcal{P}_r(n)$, then
we write $\blam\vdash n$ and $|\ulam|=n$. Then $\mathcal{P}_{r}(n)$
is a poset under the \textbf{dominance partial order} ``$\unrhd$",
where $\ulam\gedom\umu$ if
$$\sum_{a=1}^{s-1}|\lambda^{(a)}|+\sum_{j=1}^i\lambda^{(s)}_j
\ge\sum_{a=1}^{s-1}|\mu^{(a)}|+\sum_{j=1}^i\mu^{(s)}_j,$$ for all
$1\le s\le r$ and all $i\ge1$.

Let $\ulam\in\mathcal{P}_r(n)$. Recall that the \textbf{Young
diagram} of~$\ulam$ is the set $$
[\ulam]=\bigl\{(i,j,s)\bigm|\text{$1\le
j\le\lambda^{(s)}_i$}\bigr\}.$$ The elements of $[\ulam]$ are called
the nodes of $\ulam$. A \textbf{$\ulam$-tableau} is a bijection
$\t\colon {[\ulam]}\rightarrow\{1,2,\dots,n\}$. The $\ulam$--tableau
$\t$ is \textbf{standard} if $\t(i,j,s)\leq\t(i',j',s)$ whenever
$i\le i'$, $j\le j'$. Let $\std(\ulam)$ be the set of standard
$\ulam$--tableaux. For any two nodes $\gamma=(a,b,c),
\gamma'=(a',b',c')$ of $\ulam$, say that $\gamma$ is {\it below}
$\gamma'$, or $\gamma'$ is {\it above} $\gamma$, if either $c>c'$ or
$c=c'$ and $a>a'$. If $\gamma'$ is above $\gamma$ then we write $\gamma'>\gamma$. A removable node of $\ulam$ is a triple
$(i,j,s)\in[\blam]$ such that $[\blam]-\{(i,j,s)\}$ is the Young
diagram of a multipartition, while an addable node of $\ulam$ is a
triple $(i,j,s)$ which does not lie in $[\ulam]$ but is such that
$[\ulam]\cup\{(i,j,s)\}$ is the Young diagram of a multipartition.

Now let $e$ be $0$ or an integer bigger than $1$. Let $v_1,\cdots,v_r\in\Z/e\Z$.
Let $\bQ:=(v_1,\cdots,{v_r})$. The \textbf{residue} of the node $\gamma=(a,b,c)$ is
defined to be
$$
\res(\gamma):=b-a+v_c+e\mathbb{Z}\in\mathbb{Z}/e\mathbb{Z},
$$
In this case, we say that $\gamma$ is a $\res(\gamma)$-node.

If $\umu=(\mu^{(1)},\cdots,\mu^{(r)})$ is an $r$-partition of
$n+1$ with $[\umu]=[\ulam]\cup\bigl\{\gamma\bigr\}$ for some
removable node $\gamma$ of $\umu$, we write $\ulam\rightarrow\umu$
or $\bmu/\blam=\gamma$. If in addition $\res(\gamma)=x$, we also
write $\ulam\overset{x}{\rightarrow}\umu$.

\begin{dfn} {\rm (\!\cite{AM})} Let $x\in\Z/e\Z$. Let $\blam\in\mathcal{P}_r(n)$ and $\eta$ be a removable $x$-node of $\blam$. If whenever $\gamma$ is an addable $x$-node of $\blam$ which is
below $\eta$, there are more removable $x$-nodes between $\gamma$
and $\eta$ than there are addable $x$-nodes, then we call $\eta$ a
{\bf normal} $x$-node of $\blam$. The unique highest normal $x$-node
of $\blam$ is called the {\bf good} $x$-node of $\blam$;
\end{dfn}

For example, suppose $n=19, r=3, e=4$, $v_1=4\Z, v_2=2+4\Z, v_3=4\Z$.
The nodes of $\ulam=((2),(4,2,2),(5,2,1,1))$ have the following
residues
$$
\ulam=\biggl(\left(\begin{matrix} 0& 1
\end{matrix}
\right),\left(\begin{matrix} 2& 3& 0& 1\\
1& 2& &\\
0& 1& &
\end{matrix}
\right),\left(\begin{matrix} 0& 1& 2& 3& 0\\
3& 0& & &\\
2& & & &\\
1& & & &
\end{matrix}
\right)\biggr).
$$
$\blam$ has six removable nodes. Fix a residue $x$ and consider the
sequence of removable and addable $x$-nodes obtained by reading the
boundary of $\blam$ from the bottom up. In the above example, we
consider the residue $x=1$, then we get a sequence ``RAARRR", where each
``A'' corresponds to an addable $x$-node and each ``R'' corresponds
to a removable $x$-node. Given such a sequence of letters A,R, we
remove all occurrences of the string ``AR'' and keep on doing this
until no such string ``AR'' is left. The {\it normal} $x$-nodes of
$\ulam$ are those that correspond to the remaining ``R'' and the
highest of these is the {\it good} $x$-node. In the above example,
there are two normal $1$-nodes $(1,2,1)$ and $(4,1,3)$. The removable $1$-node $(1,2,1)$ is a
good $1$-node. If $\gamma$ is a good $x$-node of $\bmu$ and $\ulam$
is the multipartition such that $[\umu]=[\ulam]\cup\gamma$, we write
$\ulam\overset{x}{\twoheadrightarrow}\umu$.

\begin{dfn} {\rm(\!\!\cite{AM})} The set $\mathcal{K}_r(n)$ of Kleshchev $r$-partitions with respect to
$(e;\bQ)$ is defined inductively as follows:
\begin{enumerate}
\item[(1)] $\mathcal{K}_r(0):=\Bigl\{\uempty:=\bigl(\underbrace{\emptyset,\cdots,
\emptyset}_{\text{$r$ copies}}\bigl)\Bigr\}$;
\item[(2)] $\mathcal{K}_{r}(n+1):=\Bigl\{\umu\in\mathcal{P}_{r}(n+1)\Bigm|\text{$\ulam\overset{x}
{\twoheadrightarrow}\umu$ for some $\ulam\in\mathcal{K}_r(n)$ and
$x\in\Z/e\Z$}\Bigr\}$. \end{enumerate} {\it Kleshchev's good
lattice} with respect to $(e;\bQ)$ is the infinite graph whose
vertices are the Kleshchev $r$-partitions with respect to
$(e;\bQ)$ and whose arrows are given by $$
\text{$\ulam\overset{x}{\twoheadrightarrow}\umu$\quad$\Longleftrightarrow$\quad
$\ulam$ is obtained from $\umu$ by removing a good $x$-node}. $$
\end{dfn}

Let $K$ be a field. Let $q$ be a primitive $e$-th root of unity if
$e>1$; or not a root of unity if $e=0$. The Ariki--Koike algebra
$\HH_{r,n}(q;q^{v_1},\cdots,q^{v_r})$ (or the cyclotomic Hecke
algebra of type $G(r,1,n)$) is the associative unital $K$-algebra
with generators $T_0,T_1,\cdots,T_{n-1}$ and relations
$$\begin{aligned}
&(T_{0}-q^{v_1})\cdots (T_{0}-q^{v_r})=0,\\
&T_0T_1T_0T_1=T_1T_0T_1T_0,\\
&(T_i+1)(T_i-q)=0,\quad\text{for $1\leq i\leq n-1$,}\\
&T_iT_{i+1}T_i=T_{i+1}T_{i}T_{i+1},\quad\text{for $1\leq i\leq n-2$,}\\
&T_iT_j=T_jT_i,\quad\text{for $0\leq i<j-1\leq n-2$.}\end{aligned}
$$
These algebras were introduced in the work of Brou\'e and
Malle~\cite{BM} and of Ariki and Koike~\cite{AK}. They include the
Iwahori--Hecke algebras of types $A$ and $B$ as special cases.
Conjecturally, they have an intimate relationship with the
representation theory of finite reductive groups. The modular
representation theory of these algebras was studied in \cite[Section
5]{GL} and \cite{DJM}, where $\HH_{r,n}(q;q^{v_1},\cdots,q^{v_r})$
was shown to be a cellular algebra in the sense of \cite{GL}. Using
the cellular basis constructed in \cite{DJM}, we know that the
resulting cell modules (i.e., Specht modules)
$\{S_{\ulam}\}_{\ulam\vdash n}$ are indexed by the set of
$r$-partitions of $n$. By the theory of cellular algebras, each
Specht module $S_{\ulam}$ is equipped with a bilinear form
$\langle,\rangle$. Let $D_{\ulam}:=S_{\ulam}/\rad\langle,\rangle$.
The set $\bigl\{D_{\ulam}\bigm|D_{\ulam}\neq 0, \ulam\vdash
n\bigr\}$ is a complete set of pairwise non-isomorphic absolutely
simple $\HH_{r,n}(q;q^{v_1},\cdots,q^{v_r})$-modules. The
significance of the notion of Kleshchev multipartition can be seen
from the following remarkable result of Ariki.

\begin{thm} {\rm(\!\!\cite[Theorem 4.2]{A3})}\label{athm} Let $\blam\in\mathcal{P}_r(n)$. Then, $D_{\blam}\neq 0$ if and only
if $\blam\in\mathcal{K}_r(n)$.
\end{thm}

\begin{dfn} \label{rd} Let $\blam\in\mathcal{P}_r(n)$ and $\t\in\std(\blam)$. The residue sequence of\, $\t$ is defined
to be the sequence
$$\bigl(\res(\t^{-1}(1)),\cdots,\res(\t^{-1}(n))\bigr).$$
\end{dfn}

The following definitions are natural extensions of the
corresponding definitions given in the case where $r=1,2$, see
\cite{AJ}, \cite{DJM2} and \cite{GL}.

\begin{dfn} Let $\blam\in\mathcal{P}_r(n)$. $\blam$ is said to be $(\bQ,e)$-restricted if there exists $\t\in\std(\blam)$
such that the residue sequence of any standard tableau of shape
$\bmu\lhd\blam$ is not the same as the residue sequence of $\t$.
\end{dfn}

Note that if $r=1$, by \cite[Corollary 3.41]{Ma1},
$(\bQ,e)$-restricted partitions are the same as $e$-restricted
partitions. In particular, $\mathcal{K}_1(n)$ is the same as the set
of $e$-restricted partitions. If $r=2$, the above definition
appeared in the paper \cite{DJM2} of Dipper--James--Murphy. They
proved that if $\blam$ is $(\bQ,e)$-restricted, then $D_{\blam}\neq
0$, and they conjectured the converse is also true, i.e.,
$\blam\in\mathcal{K}_2(n)$ if and only if $\blam$ is
$(\bQ,e)$-restricted. This conjecture was recently proved by
Ariki--Jacon \cite{AJ}, using a new characterization of Kleshchev
bipartitions obtained in \cite{AKT}. The general case (i.e., when
$r>2$) remains open. That is

\begin{conj0} \label{GDJM} Let $\blam\in\mathcal{P}_r(n)$. Then $\blam\in\mathcal{K}_r(n)$ if and only if $\blam$ is
$(\bQ,e)$-restricted.
\end{conj0}

Note that the generalised DJM conjecture can be understood as a
criterion for $D_{\blam}$ to be non-zero, where $D_{\blam}$ is defined using
the Dipper--James--Mathas cellular basis of
$\HH_{r,n}(q;q^{v_1},\cdots,q^{v_r})$. With respect to a different
cellular basis, Graham and Lehrer proposed a similar conjecture in
\cite[(5.9),(5.10)]{GL}. Since we do not know whether the two set of
cellular datum give rise to equivalent cell modules and labeling of
simple modules when $r>2$, it is not clear to us whether the two
conjectures are equivalent or not.

In fact, the ``if" part of Conjecture \ref{GDJM} is easy, as we
shall describe in the following. The definition of
$(\bQ,e)$-restricted multipartition can be reformulated in terms of
the action of the affine quantum group on a Fock space (cf.
\cite{AJ}). To recall this, we need some more notations. Let $v$ be
an indeterminate over $\Q$. Let $\mathfrak{g}:=\ksl_e$ be the affine
Lie algebra of type $A_{e-1}^{(1)}$ if $e>1$; or
$\mathfrak{g}:=\mathfrak{gl}_{\infty}$ be the affine Lie algebra of
type $A_{\infty}$ if $e=0$. Let $U_v(\mathfrak{g})$ be the
corresponding affine quantum group with Chevalley generators $e_i, f_i, k_i$ and $k_d$ for $i\in\Z/e\Z$. 
Let $\bigl\{\Lambda_i\bigm|i\in\Z/e\Z\bigr\}$ be the set of fundamental
weights of $\mathfrak{g}$. Let $\mathcal{F}$ be the level $r$
$v$-deformed Fock space associated to $(e;v_1,\cdots,v_r)$ which was used in \cite{AM}.
Our space $\mathcal{F}$ was denoted by $\mathcal{F}_v$ in \cite{AM} and one should understand 
the $r$-tuple $(Q_1,\cdots,Q_r)$ in \cite{AM} as the $r$-tuple $(q^{v_1},\cdots,q^{v_r})$ in this paper, where 
$q$ is a primitive $e$-th root of unity in $\mathbb{C}$ if $e>1$; or not a root of unity if $e=0$.  
By definition\footnote{Although in \cite{AM} the ground field of $\mathcal{F}_v$ is $\mathbb{C}(v)$, it does no harm to replace it by $\Q(v)$.}, $\mathcal{F}$ is a
$\Q(v)$-vector space with the basis given by the set of all
$r$-partitions, i.e., $$ \mathcal{F}=\bigoplus_{n\geq 0,
\blam\in\mathcal{P}_r(n)}\Q(v)\blam.
$$
By \cite{MM} and \cite{AM}, there is an action of $U_v(\mathfrak{g})$ on $\mathcal{F}$ which
quantizes the classical action of $\mathfrak{g}$ on the $\Q$-vector
space $\bigoplus_{n\geq 0, \blam\in\mathcal{P}_r(n)}\Q\blam$. That is, for each $i\in\Z/e\Z$ and $\blam\in\mathcal{P}_r(n)$,
$$\begin{aligned}
&e_i\blam=\sum_{\bmu\overset{i}\rightarrow\blam}v^{-N_i^{r}(\bmu,\blam)}\bmu,\quad
f_i\blam=\sum_{\blam\overset{i}\rightarrow\bmu}v^{N_i^{l}(\blam,\bmu)}\bmu,\\
& k_i\blam=v^{N_i(\blam)}\blam,\quad k_d\blam=v^{-N_d(\blam)}\blam,
\end{aligned}
$$
where $$\begin{aligned}
&N_i^{r}(\bmu,\blam):=\#\Biggl\{\gamma\Biggm|\begin{matrix}\text{$\gamma$ is an addable $i$-node}\\
\text{for $\blam$, $\gamma>\blam/\bmu$}\end{matrix}\Biggr\}-\#\Biggl\{\gamma\Biggm|\begin{matrix}\text{$\gamma$ is a removable $i$-node}\\
\text{for $\blam$, $\gamma>\blam/\bmu$}\end{matrix}\Biggr\},\\
&N_i^{l}(\blam,\bmu):=\#\Biggl\{\gamma\Biggm|\begin{matrix}\text{$\gamma$ is an addable $i$-node}\\
\text{for $\blam$, $\gamma<\bmu/\blam$}\end{matrix}\Biggr\}-\#\Biggl\{\gamma\Biggm|\begin{matrix}\text{$\gamma$ is a removable $i$-node}\\
\text{for $\blam$, $\gamma<\bmu/\blam$}\end{matrix}\Biggr\},\\
&N_i(\blam)=\#\Bigl\{\gamma\Bigm|\begin{matrix}\text{$\gamma$ is an addable}\\
\text{ $i$-node for $\blam$}\end{matrix}\Bigr\}-\#\Bigl\{\gamma\Bigm|\begin{matrix}\text{$\gamma$ is a removable}\\
\text{ $i$-node for $\blam$}\end{matrix}\Bigr\},\\
&N_d(\blam):=\#\Bigl\{\gamma\in[\blam]\Bigm|\res(\gamma)=0\Bigr\}.
\end{aligned}
$$
Note that the empty multipartition $\uempty$ is a highest weight vector of weight $\sum_{j=1}^{r}\Lambda_{v_{j}}$. 
One can also identify $\mathcal{F}$ with a tensor product of $r$ level one Fock spaces. We refer the reader to the proof of \cite[Proposition 2.6]{AM} for more details.

\begin{lem} \label{2dfn} Let $\blam\in\mathcal{P}_r(n)$. Then $\blam$ is $(\bQ,e)$-restricted if and only if there exists a sequence
$(i_1,\cdots,i_n)$ of residues such that $$
f_{i_n}\cdots f_{i_1}\uempty=A_0\blam+\sum_{\bmu\ntriangleleft\blam}A_{\blam,\bmu}(v)\bmu,
$$
for some $A_0, A_{\blam,\bmu}(v)\in\Z_{\geq 0}[v,v^{-1}]$ with $A_0\neq 0$, where
$f_{i_1},\cdots, f_{i_n}$ are the Chevalley generators of
$U_v(\mathfrak{g})$.
\end{lem}
\begin{proof} For any residue $j$ and any $\bmu\in\mathcal{P}_r(n)$, by definition, we have that
$$
f_j\bmu=\sum_{\res(\bnu/\bmu)=j}C_{\bmu,\bnu}(v)\nu,
$$
for some $C_{\bmu,\bnu}(v)\in\Z_{\geq 0}[v,v^{-1}]$ satisfying
$C_{\bmu,\bnu}(1)\neq 0$. The lemma follows directly from this fact and
the definition of standard tableaux.
\end{proof}

\begin{cor} \label{mcor1} Let $\blam\in\mathcal{P}_r(n)$. If $\blam$ is $(\bQ,e)$-restricted, then
$\blam\in\mathcal{K}_r(n)$.
\end{cor}

\begin{proof} Let $\Lambda:=\sum_{j=1}^{r}\Lambda_{v_{j}}$. Then the $U_v(\mathfrak{g})$-submodule
of $\mathcal{F}$ generated by $\uempty$ is the
irreducible highest weight module $V(\Lam)$ of highest weight
$\Lam$. It is well-known that $V(\Lam)$ has a canonical basis
$\bigl\{G(\mu)\bigr\}$, which is indexed by the set $\mathcal{K}_r:=
\bigsqcup_{n\geq 0}\mathcal{K}_r(n)$. Combining \cite[Theorem 4.4]{A2} with \cite[Theorem 3.14, (5.3), Theorem 5.6, Theorem 5.14, Corollary 5.15]{BK}, we see that for any
$\bmu\in\mathcal{K}_r(n)$, 
$$
G(\bmu)=\bmu+\sum\limits_{\substack{\bnu\in\mathcal{P}_r(n)\\
\bnu\rhd\bmu}}d_{\bnu,\bmu}(v)\bnu,$$ where $d_{\bnu,\bmu}(v)\in\mathbb{Z}_{\geq 0}[v,v^{-1}]$ for each 
$\bnu\in\mathcal{P}_r(n)$ satisfying $\bnu\rhd\bmu$.

Since $\blam$ is $(\bQ,e)$-restricted, by Lemma \ref{2dfn}, we deduce that there exists a sequence
$(i_1,\cdots,i_n)$ of residues such that $$
f_{i_n}\cdots f_{i_1}\uempty=A_0\blam+\sum_{\bmu\ntriangleleft\blam}A_{\blam,\bmu}(v)\bmu,
$$
for some $A_0, A_{\blam,\bmu}(v)\in\Z_{\geq 0}[v,v^{-1}]$ with $A_0\neq 0$.

On the other hand, since $f_{i_n}\cdots f_{i_1}\uempty\in L(\Lam)$, hence we can write $$
A_0\blam+\sum_{\bmu\ntriangleleft\blam}A_{\blam,\bmu}(v)\bmu=\sum_{\nu\in\mathcal{K}_r(n)}A'_{\nu}(v)G(\nu),
$$
for some $A'_{\nu}(v)\in\Z[v,v^{-1}]$. It follows from the induction
on the dominance partial order ``$\lhd$" that
$\blam\in\mathcal{K}_r(n)$, as required.
\end{proof}

Therefore, we have proved that ``if" part of the Conjecture
\ref{GDJM}. It remains to consider the ``only if" part of that
conjecture. To this end, we need a result of Littelmann.

We need some more notations. Let
$P^{+}:=\bigl\{\sum_{i\in\Z/e\Z}a_i\Lambda_i\bigm|a_i\in\Z_{\geq
0},\forall\,i\bigr\}$ be the set of dominant weights. Let
$\{\alpha_i\}_{i\in\Z/e\Z}$ (resp., $\{h_i\}_{i\in\Z/e\Z}$) be the
set of simple roots (resp., simple coroots). For each dominant
weight $\Lambda$, let $V(\Lambda)$ be the irreducible
$U_v(\mathfrak{g})$-module with highest weight $\Lambda$. We assume
that the reader is familiar with the theory of Kashiwara crystals.
It is well-known that $V(\Lambda)$ has a crystal basis. We denote by
$B(\Lambda)$ its crystal graph. Note that $B(\Lambda)$ is equipped
with additional data $\wt$, $\eps_i$, $\varphi_i$, $\widetilde{e}_i$
and $\widetilde{f}_i$. We refer the readers to \cite{A4}, \cite{Ka0} and \cite{Ka} for
details. We use $u_{\Lambda}$ to denote the unique element in
$B(\Lambda)$ satisfying $\wt(u_{\Lambda})=\Lambda$. For each
$i\in\Z/e\Z$, there are two important realizations of the crystal
graph $B(\Lambda_{i})$, one by $e$-restricted partitions (cf.
\cite{MM}), the other by Littelmann's path model (cf. \cite{L1}). Let $W$ be the
affine Weyl group with standard Coxeter generators $s_i$,
$i\in\Z/e\Z$. By definition, $W$ is presented by the generators
$s_i$, $i\in\Z/e\Z$ and the following relations:
$$\begin{aligned}
&s_i^2=1, \quad\forall\, i\in\Z/e\Z;\\
&s_is_j=s_js_i,\quad\text{if $i\neq j\pm 1$;}\\
&s_is_{i+1}s_i=s_{i+1}s_is_{i+1}, \quad\forall\, i\in\Z/e\Z.
\end{aligned}
$$
With these two realizations in mind, we can associate each
$e$-restricted partition $\lam$ with an Lakshimibai-Seshadri path
$(w_1\Lambda_i,\cdots,{w}_s\Lambda_i; a_0,\cdots,a_s)$, where
$w_1,\cdots,w_s$ are elements in $W$ such that
$w_1\Lambda_i,\cdots,{w}_s\Lambda_i$ are distinct and
$0=a_0<a_1<\cdots<a_s=1$ are some rational numbers. We
refer the readers to \cite{L1} for the precise definition of
Lakshimibai--Seshadri paths and related notions. For simplicity, we shall often abbreviate ``Lakshimibai--Seshadri paths" to ``LS paths". Note that there
is a canonical crystal embedding
$B(\Lambda_{v_1}+\cdots+\Lambda_{v_r})\hookrightarrow
B(\Lambda_{v_1})\otimes\cdots\otimes B(\Lambda_{v_r})$. We identify
$B(\Lambda_{v_1}+\cdots+\Lambda_{v_r})$ with the image by this
embedding. We have the following result of Littelmann, which was
reformulated in \cite[Theorem 5.7]{AKT}.

\begin{thm} {\rm(\!\!\cite[Theorem 10.1]{L2})} \label{pthm} Let $$
\pi=\pi^{(1)}\otimes\cdots\otimes\pi^{(r)}\in
B(\Lambda_{v_1})\otimes\cdots\otimes B(\Lambda_{v_r}).
$$
Then $\pi$ belongs to $B(\Lambda_{v_1}+\cdots+\Lambda_{v_r})$ if and
only if there exists a sequence $$ w_{1}^{(1)}\geq\cdots\geq
w_{N_1}^{(1)}\geq w_{1}^{(2)}\geq\cdots\geq
w_{N_2}^{(2)}\geq\cdots\geq w_{N_r}^{(r)}
$$
in $W$ such that $$
\pi^{(k)}=(w_1^{(k)}\Lambda_{v_k},\cdots,w_{N_k}^{(k)}\Lambda_{v_k};a_0^{(k)},\cdots,a_{N_k}^{(k)}),
$$
for any integer $1\leq k\leq r$, where $``\geq"$ is the Bruhat order
and $0=a_0^{(k)}<a_1^{(k)}<\cdots<a_{N_k}^{(k)}=1$ are some rational
numbers.
\end{thm}

\begin{lem} {\rm (\!\!\cite[Proposition 4.10]{Ma2})} \label{eres} If $\ulam=(\lam^{(1)},\cdots,\lam^{(r)})$ is a Kleshchev $r$-partition with respect to
$(e;{v_1},\cdots,{v_r})$, then for each $1\leq i\leq r$,
$\lam^{(i)}$ is an $e$-restricted partition.
\end{lem}

The combinatorial realization of $B(\Lambda_i)$ in terms of
$e$-restricted partitions (cf. \cite{MM}) allows a natural
generalization to the higher level case (cf. \cite{AM}) as we now
recall. Let $\Lambda:=\Lambda_{v_1}+\cdots+\Lambda_{v_r}$. Set
$\mathcal{K}_r:=\sqcup_{n\geq 0}\mathcal{K}_r(n)$. For each
$i\in\Z/e\Z$ and $\blam\in\mathcal{K}_r(n)$, we define
$$\begin{aligned}
\widetilde{e}_i\blam&=\begin{cases}\blam-\{\gamma\}, &\text{if
$\blam$ has a good
$i$-node $\gamma$;}\\
0, &\text{otherwise;}
\end{cases}\\
\widetilde{f}_i\blam&=\begin{cases}\blam\cup\{\gamma\}, &\text{if
$\gamma$ is a good $i$-node of $\blam\cup\{\gamma\}$;}\\
0, &\text{otherwise.}
\end{cases}\\
\varepsilon_i(\blam)&=\max\{n\geq 0|\widetilde{e}_i^nb\neq
0\},\,\,\, \varphi_i(\blam)=\max\{n\geq 0|\widetilde{f}_i^nb\neq
0\};\\
\wt(\blam)&=\Lambda-\sum_{i\in\Z/e\Z}\widetilde{N}_i(\blam)\alpha_i,
\end{aligned}
$$
where $\widetilde{N}_i(\blam)$ is the number of $i$-nodes in $[\blam]$. By a
result of Misra--Miwa \cite{MM} and Ariki--Mathas \cite{AM}, the
data $\mathcal{K}_r$, $\wt$, $\eps_i$, $\varphi_i$,
$\widetilde{e}_i$ and $\widetilde{f}_i$ define a realization of the
crystal $B(\Lambda)$ in terms of Kleshchev's good lattice with
respect to $(e,\bQ)$. Henceforth, we make this identification. In
particular, taking $r=1$ and $i\in\Z/e\Z$, we can identify any
element in $B(\Lambda_i)$ with an $e$-restricted partition.

Recall that for any two $\mathfrak{g}$-crystals $B_1, B_2$. The
tensor product $B_1\otimes B_2$ is the set $B_1\times B_2$ equipped
with the crystal structure defined by \begin{enumerate}
\item $\wt(b_1\otimes b_2)=\wt(b_1)+\wt(b_2)$;
\item $$
\widetilde{e}_i(b_1\otimes b_2)=\begin{cases}
\widetilde{e}_ib_1\otimes b_2, &\text{if
$\varphi_i(b_1)\geq\epsilon_i(b_2)$;}\\
b_1\otimes \widetilde{e}_ib_2, &\text{if
$\varphi_i(b_1)<\epsilon_i(b_2)$;}\\
\end{cases}
$$
\item $$
\widetilde{f}_i(b_1\otimes b_2)=\begin{cases}
\widetilde{f}_ib_1\otimes b_2, &\text{if
$\varphi_i(b_1)>\epsilon_i(b_2)$;}\\
b_1\otimes \widetilde{f}_ib_2, &\text{if
$\varphi_i(b_1)\leq\epsilon_i(b_2)$.}\\
\end{cases}
$$
\item $\epsilon_i(b_1\otimes b_2)=\max\bigl\{\epsilon_i(b_1),\,\,\epsilon_i(b_2)-\langle h_i,\wt(b_1)\rangle\bigr\}$;
\item $\varphi_i(b_1\otimes b_2)=\max\bigl\{\varphi_i(b_1)+\langle
h_i,\wt(b_2)\rangle,\,\,\varphi_i(b_2)\bigr\}$.
\end{enumerate}

Let $\blam=(\lam^{(1)},\cdots,\lam^{(r)})\in\mathcal{P}_r(n)$. By
Lemma \ref{eres}, $\blam$ is a Kleshchev multipartition with respect
to $(e;\bQ)$ only if each $\lam^{(s)}$ is an $e$-restricted
partition for $1\leq s\leq r$.

\begin{lem} The map
which sends each
$\blam=(\lam^{(1)},\cdots,\lam^{(r)})\in\mathcal{K}_r(n)$ to
$\lam^{(r)}\otimes\lam^{(r-1)}\otimes\cdots\otimes\lam^{(1)}$
coincides with  the canonical crystal embedding
$B(\Lambda)=B(\Lambda_{v_r}+\cdots+\Lambda_{v_1})\hookrightarrow
B(\Lambda_{v_r})\otimes\cdots\otimes B(\Lambda_{v_1})$. In
particular, if each $\lam^{(s)}$ is $e$-restricted for $1\leq s\leq
r$, then $\ulam=(\lam^{(1)},\cdots,\lam^{(r)})$ is a Kleshchev
$r$-partition with respect to $(e;\bQ)$ if and only if
$\lam^{(r)}\otimes\cdots\otimes\lam^{(1)}$ belongs to
$B(\Lambda_{v_r}+\cdots+\Lambda_{v_1})$.
\end{lem}

\begin{proof} This follows directly from the fact that if we reverse the order of the components,
then the action of the Kashiwara operator $\widetilde{f}_i$ on
tensor product of crystals coincides with the operator of adding
good $i$-node on Kleshchev multipartitions.
\end{proof}

The above lemma implies that the problem of characterizing Kleshchev
multipartition in terms of its components is essentially a purely
crystal theoric question. For the latter, Theorem \ref{pthm}
gives an answer in the language of Littelmann's path model. The
following lemma is a direct consequence of Theorem \ref{pthm}.

\begin{cor} \label{pcor1} If $\ulam=(\lam^{(1)},\cdots,\lam^{(r)})$ is a Kleshchev $r$-partition with respect to
$(e;{v_1},\cdots,{v_r})$, then $(\lam^{(j_1)},\cdots,\lam^{(j_t)})$
is a Kleshchev $t$-partition with respect to
$(e;{v_{j_1}},\cdots,{v_{j_t}})$ for any integers $1\leq t< r$ and
$1\leq j_1<\cdots<j_t\leq r$.\end{cor}

\noindent {2.13. \it Remark.} We note that the converse of Corollary
\ref{pcor1} is in general false. For example, let $e=5,
(v_1,v_2,v_3)=(3+5\Z,2+5\Z,1+5\Z)$. Let $$ \lam^{(1)}:=(5,1),\quad
\lam^{(2)}:=(3),\quad \lam^{(3)}:=(2).
$$
Then it is easy to check that \begin{enumerate}
\item[(1)] $(\lam^{(1)},\lam^{(2)})$ is a Kleshchev bipartition with respect to $(5;v_1,v_2)$; and
\item[(2)] $(\lam^{(2)},\lam^{(3)})$ is a Kleshchev bipartition with respect to $(5;v_2,v_3)$; and
\item[(3)] $(\lam^{(1)},\lam^{(3)})$ is a Kleshchev bipartition with respect to $(5;v_1,v_3)$;
\end{enumerate}
but $(\lam^{(1)},\lam^{(2)},\lam^{(3)})$ is not a Kleshchev
$3$-partition with respect to $(5;v_1,v_2,v_3)$.
\smallskip

Let $\Lambda$ be any dominant weight. By \cite[Corollaire
8.1.5]{Ka}, there exists a unique crystal morphism $K_h:
B(\Lambda)\hookrightarrow B(h\Lambda)$ of amplitude $h$, for all
$h\in\mathbb{N}$. In other words,
\begin{enumerate}
\item[(i)] $K_{h}(u_{\Lambda})=u_{h\Lambda}$;
\item[(ii)] $\wt(K_h(b))=h\wt(b), \eps_i(K_h(b))=h\eps_i(b)$ and $\varphi_i(K_h(b))=h\varphi_i(b)$;
\item[(iii)] $K_h(\widetilde{e}_ib)=\widetilde{e}_i^{h}K_h(b)$ and $K_h(\widetilde{f}_ib)=\widetilde{f}_i^{h}K_h(b)$
for all $b\in B(\Lambda)$.
\end{enumerate}
Composing $K_h$ with the natural embedding
$B(h\Lambda)\hookrightarrow B(\Lambda)^{\otimes h}$, we get a
crystal morphism $S_h: B(\Lambda)\hookrightarrow B(\Lambda)^{\otimes
h}$ of amplitude $h$.\smallskip

Recall that, for each $w\in W$, the weight space $V(\Lam)_{w\Lambda}$ is one-dimensional. We use $u_{w\Lambda}$ to denote the unique element in $B(\Lam)$ satisfying $\wt(u_{w\Lambda})=w\Lam$.

\addtocounter{equation}{1}
\begin{lem} {\rm(\!\!\cite[Proposition 8.3.2]{Ka})} \label{keylem} Let $b\in B(\Lambda)$.
Then there is an integer $s>0$, and rational numbers
$0=a_0<a_1<\cdots<a_s=1$ and elements $w_1,\cdots,w_s$ of $W$ such
that $w_1\Lambda,\cdots,w_s\Lambda$ are pairwise distinct and
whenever $h$ satisfies $(a_{i+1}-a_{i})h\in\mathbb{Z}_{\geq 0}$ for
all $i$ we have
$$S_h(b)=u_{w_1\Lambda}^{\otimes (a_1-a_0)h}\otimes
u_{w_2\Lambda}^{\otimes (a_2-a_1)h}\otimes\cdots\otimes
u_{w_s\Lambda}^{\otimes (a_s-a_{s-1})h}.$$ Furthermore, the map
$$ b\mapsto (w_1\Lambda,\cdots,w_s\Lambda;a_0,a_1,\cdots,a_s)
$$
coincides with Littelmann's path model.
\end{lem}

{\it Henceforth, we assume that $\Lam=\Lam_k$ is a fundamental
weight.} We use $W_k$ to denote the symmetric group generated by
$s_i, i\in\Z/e\Z-\{k+e\Z\}$. Recall that the crystal $B(\Lam_k)$ has
a realization in terms of the set of $e$-restricted partitions. We
denote by $\emptyset_k$ the empty partition in $B(\Lambda_k)$. Let
$W/W_k$ be the set of distinguished coset representatives of $W_k$
in $W$. For any $i\in\Z/e\Z$ and any $e$-core $\nu$, let $s_i\nu$ be
defined as in \cite[Lemma 3.3(2)]{AKT}. For any $w\in W$,
$w\nu:=s_{i_1}\cdots s_{i_t}\nu$ if $s_{i_1}\cdots s_{i_t}$ is a
reduced expression of $w$. This is well-defined, i.e., independent
of the choice of the reduced expression. By \cite[Proposition
3.5]{AKT}, if $w\in W/W_k$, then $w\emptyset_k$ is an $e$-core, and
this gives rise to a natural bijection between the set of $e$-cores
and the set $\bigl\{w\Lam\bigm|w\in W\bigr\}$. Note that the empty
partition $\emptyset_k$ corresponds to $u_{\Lam}$, while the
$e$-core $w\emptyset$ corresponds to $u_{w\Lam}$. Therefore,
translating into the language of $e$-cores, we can write an LS-path
$(w_1\Lam,\cdots,w_s\Lam;a_0,a_1,\cdots,a_s)$ as
$(\nu_1,\cdots,\nu_s;a_0,a_1,\cdots,a_s)$, where
$\nu_i=w_i\emptyset_k$ is an $e$-core for each $i$, and we can
rephrase Lemma \ref{keylem} as follows: for each $e$-restricted
partition $\lam$, there exist an integer $s$, distinct $e$-cores
$\nu_1,\cdots,\nu_s$ and rational numbers $0=a_0<a_1<\cdots<a_s=1$
such that whenever $h$ satisfies
$(a_{i+1}-a_{i})h\in\mathbb{Z}_{\geq 0}$ for all $i$ we have
$$S_h(\lam)=\nu_1^{\otimes (a_1-a_0)h}\otimes
\nu_2^{\otimes (a_2-a_1)h}\otimes\cdots\otimes \nu_s^{\otimes
(a_s-a_{s-1})h}.$$ Furthermore, the map $\lam\mapsto
(\nu_1,\cdots,\nu_s;a_0,a_1,\cdots,a_s)$ coincides with Littelmann's
path model.

In the remaining part of this paper, we shall always write an
LS-path
$$(w_1\Lambda_k,\cdots,w_s\Lam_k;a_0,a_1,\cdots,a_s)$$ as
$(\nu_1,\cdots,\nu_s;a_0,a_1,\cdots,a_s)$, where $\nu_i$ is the
unique $e$-core such that $\nu_i=w_i\emptyset$ for each $1\leq i\leq
s$.

\begin{lem} {\rm(\!\!\cite{AKT}, \cite{Ka})} \label{keylm} With the notations as above, the map which sends each
$\lam$ to $(\nu_1,\cdots,\nu_s;a_0,a_1,\cdots,a_s)$ defines an
isomorphism of crystals between the two realizations of $B(\Lam_k)$,
the one by $e$-restricted partitions and the one by LS-paths.
Furthermore, if $\lam$ is an $e$-core, then $s=1$ and $\nu_1=\lam$.
\end{lem}
\begin{proof} The first part of the lemma follows from \cite[Theorem 5.14]{AKT} and \cite[Theorem 8.2.3]{Ka}, while 
the second part of the lemma is a direct consequence of \cite[Proposition 8.3.2 (1)]{Ka}.
\end{proof}

In the above lemma, whenever $\lam$ is mapped to $(\nu_1,\cdots,\nu_s;a_0,a_1,\cdots,a_s)$,
we then write $$\pi(\lam)=(\nu_1,\nu_2,\cdots,\nu_s).$$

\begin{lem} \label{plem} With the notations as above, we have that $\nu_1\supset\nu_2\supset\cdots\supset\nu_s$.
In particular, $(\nu_s,\cdots,\nu_1)$ is a Kleshchev
$s$-partition with respect to $(e;k,\cdots,k)$.
\end{lem}
\begin{proof} Let $\Lam=\Lam_k$. By Lemma \ref{keylem}, there exist integers $n_1,\cdots, n_s$, such that $$
\nu_{1}^{\otimes n_1}\otimes \nu_{2}^{\otimes
n_2}\otimes\cdots\otimes \nu_{s}^{\otimes n_s}\in B(m\Lam)\subset
B(\Lam)^{\otimes m},
$$
where $m=n_1+\cdots+n_s$. For each integer $1\leq j\leq s$, we write
$\nu_j=d_jW_k$ for a unique $d_j\in W/W_k$. Applying Theorem
\ref{pthm}, we deduce that $d_1\geq d_2\geq\cdots\geq d_s$. Finally,
applying \cite[Proposition 4.4]{AKT} and Lemma \ref{keylm}, we get
that $$ \nu_1\supset\nu_2\supset\cdots\supset\nu_s,
$$
and $(\nu_s,\cdots,\nu_1)$ is a Kleshchev $s$-partition with
respect to $(e;k,\cdots,k)$.
\end{proof}

Let $\ulam=(\lam^{(1)},\cdots,\lam^{(r)})$ be an $r$-partition.
Suppose that for each integer $1\leq i\leq r$, the component $\lam^{(i)}$ is an $e$-restricted partition.  For each integer
$1\leq i\leq r$, we identify $\lam^{(i)}$ as an element in $B(\Lambda_{v_i})$ and we write
$\pi(\lam^{(i)})=(\nu_1^{(i)},\cdots,\nu_{s(i)}^{(i)})$ for some
integer $s(i)$ and some pairwise distinct $e$-cores
$\nu_1^{(i)},\cdots,\nu_{s(i)}^{(i)}$. We define $$
\widetilde{\pi}(\ulam)=\bigl(\nu_1^{(1)}\otimes\cdots\otimes\nu_{s(1)}^{(1)}\bigr)\otimes\cdots\otimes\bigl(\nu_1^{(r)}\otimes\cdots\otimes\nu_{s(r)}^{(r)}\bigr).
$$
We identify $B(\Lam_{v_1}+\cdots+\Lam_{v_r})$ with its image by the natural embedding $B(\Lam_{v_1}+\cdots+\Lam_{v_r})\hookrightarrow B(\Lam_{v_1})\otimes\cdots\otimes B(\Lam_{v_r})$ and $B\bigl(s(1)\Lam_{v_1}+\cdots+s(r)\Lam_{v_r}\bigr)$ with its image by the natural embedding $B\bigl(s(1)\Lam_{v_1}+\cdots+s(r)\Lam_{v_r}\bigr)\hookrightarrow
B(\Lam_{v_1})^{\otimes s(1)}\otimes\cdots\otimes B(\Lam_{v_r})^{\otimes s(r)}$.

\begin{cor} \label{pcor2} With the notations as above, we have that $$\lam^{(1)}\otimes\cdots\otimes\lam^{(r)}\in B(\Lam_{v_1}+\cdots+\Lam_{v_r})$$
if and only if $$ \widetilde{\pi}(\ulam)\in
B\bigl(s(1)\Lam_{v_1}+\cdots+s(r)\Lam_{v_r}\bigr).$$
\end{cor}
\begin{proof} For each $1\leq i\leq r$, we identify $\lam^{(i)}$ with an LS-path $$
\pi^{(i)}=\bigl(w_1^{(i)}\Lambda_{v_i},\cdots,w_{s(i)}^{(i)}\Lambda_{v_i};a_0^{(i)},\cdots,a_{s(i)}^{(i)}\bigr).
$$

If $\lam^{(1)}\otimes\cdots\otimes\lam^{(r)}\in B(\Lam_{v_1}+\cdots+\Lam_{v_r})$, then Theorem \ref{pthm} implies that we can choose those elements $w_1^{(i)},\cdots,w_{s(i)}^{(i)}$, $1\leq i\leq r$, in a way such that $$ w_{1}^{(1)}\geq\cdots\geq
w_{s(1)}^{(1)}\geq w_{1}^{(2)}\geq\cdots\geq
w_{s(2)}^{(2)}\geq\cdots\geq w_{s(r)}^{(r)}.
$$
By Lemma \ref{keylm}, we know that for each $1\leq i\leq r$ and $1\leq j\leq s(i)$, $\bigl(w_j^{(i)};0,1\bigr)$ is an LS-path for the $e$-core $\nu_j^{(i)}$. Applying Theorem \ref{pthm} again, we prove that $\widetilde{\pi}(\ulam)\in
B\bigl(s(1)\Lam_{v_1}+\cdots+s(r)\Lam_{v_r}\bigr)$.

Conversely, assume that $\widetilde{\pi}(\ulam)\in
B\bigl(s(1)\Lam_{v_1}+\cdots+s(r)\Lam_{v_r}\bigr)$. Then Theorem \ref{pthm} implies that we can find elements $\hat{w}_1^{(i)},\cdots,\hat{w}_{s(i)}^{(i)}$, $1\leq i\leq r$ such that $\hat{w}_j^{(i)}\Lambda_{v_i}=w_j^{(i)}\Lambda_{v_i}$ for each $1\leq i\leq r$, $1\leq j\leq s(i)$, and $$ \hat{w}_{1}^{(1)}\geq\cdots\geq
\hat{w}_{s(1)}^{(1)}\geq \hat{w}_{1}^{(2)}\geq\cdots\geq
\hat{w}_{s(2)}^{(2)}\geq\cdots\geq \hat{w}_{s(r)}^{(r)}.
$$
Since $\hat{w}_j^{(i)}\Lambda_{v_i}=w_j^{(i)}\Lambda_{v_i}$, we have that $$
\pi^{(i)}=\bigl(\hat{w}_1^{(i)}\Lambda_{v_i},\cdots,\hat{w}_{s(i)}^{(i)}\Lambda_{v_i};a_0^{(i)},\cdots,a_{s(i)}^{(i)}\bigr).
$$
Applying Theorem \ref{pthm} again, we deduce that $\lam^{(1)}\otimes\cdots\otimes\lam^{(r)}\in B(\Lam_{v_1}+\cdots+\Lam_{v_r})$ as required.
\end{proof}

\begin{cor} \label{pcor3} With the notations as above, we have that $(\lam^{(r)},\cdots,\lam^{(1)})$ is a Kleshchev
$r$-partition with respect to $(e;{v_r},\cdots,{v_1})$ if and
only if $$
\Bigl(\nu_{s(r)}^{(r)},\cdots,\nu_1^{(r)},\cdots,\nu_{s(1)}^{(1)},\cdots,\nu_1^{(1)}\Bigr)
$$
is a Kleshchev $(\sum_{i=1}^{r}s(i))$-partition with respect to
$$ \Bigl(e;\underbrace{{v_r},\cdots,{v_r}}_{\text{$s(r)$
copies}},\cdots,\underbrace{{v_1},\cdots,{v_1}}_{ \text{$s(1)$
copies}}\Bigr).
$$
\end{cor}

\bigskip
\section{The multi-core case}

The purpose of this section is to give a proof of the ``only if"
part of the Conjecture \ref{GDJM} in the multi-core case. In
particular, we prove that the generalised Dipper--James--Murphy
conjecture is true if $e=0$.\footnote{Note that in the case $e=0$,
i.e., $q$ is not a root of unity, the Ariki--Koike algebra
$\HH_{r,n}(q;q^{v_1},\cdots,q^{v_r})$ is NOT necessarily semisimple
whenever $r\geq 2$. We refer the reader to \cite{A1} for the semisimplicity criterion of Ariki--Koike algebra.}
\smallskip

\begin{dfn} \label{laddfn} Let $\blam\in\mathcal{P}_r(n)$ and
$\gamma\in[\blam]$. We call $\gamma$ a semi-ladder node of $\blam$
if $\gamma$ is a removable node of $\blam$ and there is no lower
addable node of the same residue. We call $\gamma$ a ladder node of
$\blam$ if $\gamma$ is a semi-ladder node such that there is no
higher semi-ladder node of the same residue.
\end{dfn}

A semi-ladder node of $\ulam$ is necessarily a normal node of
$\ulam$. In general, a multipartition may have no semi-ladder nodes.
Let $x\in\Z/e\Z$. By an $x$-sequence of $\blam$ we mean a sequence of removable $x$-nodes 
of $\blam$, arraying in decreasing order, i.e., $\alpha_1>\alpha_2>\cdots>\alpha_s$. 
If $\ulam$ has a semi-ladder $x$-node for some $x\in\Z/e\Z$, then we
call the sequence of all the semi-ladder $x$-nodes, arraying in
decreasing order, i.e., $\alpha_1>\alpha_2>\cdots>\alpha_s$, as a
ladder $x$-sequence of $\ulam$. It is readily seen that every node
in the ladder $x$-sequence of $\ulam$ is necessarily a normal
$x$-node of $\ulam$.

\begin{lem} \label{1ladlem} Assume that $\ulam$ is a non-empty Kleshchev
multipartition with respect to $(e;{v_1},\cdots,{v_r})$. Then
$\ulam$ has at least one ladder node.
\end{lem}
\begin{proof} By assumption, $\ulam=(\lam^{(1)},\cdots,\lam^{(r)})\in\mathcal{K}_r(n)$.
Let $\alpha$ be the lowest removable node of $\ulam$. Suppose that
$\alpha$ is in the $a$th row of the $c$th component of $\ulam$. Let
$x:=\res(\alpha)\in\Z/e\Z$. Then $\lam^{(t)}=\emptyset$ for any
integer $t>c$, and $\lam^{(c)}_s=0$ for any integer $s>a$.

By Lemma \ref{eres}, $\ulam\in\mathcal{K}_r(n)$ implies that each
component $\lam^{(t)}$ is an $e$-restricted partition. In
particular, $\lam^{(c)}$ is a non-empty $e$-restricted partition. It
follows that the residue of the unique addable node below $\alpha$
in $\lam^{(c)}$ is different from $x$. In other words, there are no
addable $x$-node in $\lam^{(c)}$ which is below $\alpha$.

We claim that for any integer $t>c$ we must have that $x\neq {v_t}$.
Since $\blam\in\mathcal{K}_r(n)$, by definition, we can find
$\blam(i)\in\mathcal{K}_r(i)$ for each $0\leq i\leq n$ such that
\begin{enumerate}
\item $\blam(0)=\uempty$, $\blam(n)=\blam$;
\item $\blam(i)\subset\blam(i+1)$ for each $0\leq i\leq n-1$;
\item $\gamma(i):=\blam(i+1)/\blam(i)$ is a good node of
$\blam(i+1)$.
\end{enumerate}
Let $1\leq j\leq n$ be the unique integer such that
$\alpha=\gamma(j)$. Since $\alpha$ is the lowest removable node of
$\ulam$, $\alpha$ must also be the lowest removable node of
$\ulam(i+1)$, it follows from the definition of good nodes that
$x\neq {v_t}$ for any $x>t$. This proves our claim.

Therefore there are actually no addable $x$-nodes in $[\ulam]$ that
can be lower than $\alpha$. Hence $\alpha$ is a semi-ladder $x$-node
and $\blam$ must have a ladder $x$-node.
\end{proof}

\begin{example} Suppose that $e=3$ and $v_1=2+3\Z, v_2=1+3\Z, v_3=3\Z$. Let $$
\blam:=\Bigl((1),\,\,(2,1),\,\,(3,1)\Bigr),\quad\,\bmu:=\Bigl((1),\,\,(3,1),\,\,\emptyset\Bigr).
$$ 
Then it is easy to see that $\blam$ is a Kleshchev $3$-partition with respect to $(3; v_1,v_2,v_3)$, while 
$\bmu$ is not a Kleshchev $3$-partition with respect to $(3; v_1,v_2,v_3)$. Furthermore, $\blam$ has a unique good node, that is, the good $2$-node $(1,1,1)$ and $\blam$ has only semi-ladder $2$-nodes, $(1,3,3)>(2,1,3)$ is the ladder $2$-sequence of $\blam$ and $(1,3,3)$ is the ladder $2$-node of $\blam$. Note that $(1,1,1)>(1,3,3)>(2,1,3)$
are all the normal $2$-nodes of $\blam$. Finally, $\bmu$ has no semi-ladder nodes.
\end{example}

\begin{dfn} \label{ladm0} Let $\blam\in\mathcal{P}_r(n)$. We call $\blam$ a ladder
multipartition if there is a sequence of residues $i_1,\cdots,i_n$
 and a sequence of multipartitions
$\blam[0]=\uempty,\cdots,\blam[n]=\blam$ such that for
each $t$, $\blam[t-1]$ is obtained from $\blam[t]$ by removing the
ladder $i_t$-node of $\blam[t]$.
\end{dfn}

\begin{dfn} \label{ladm} Let $\blam\in\mathcal{P}_r(n)$. We call $\blam$ a strong ladder
multipartition if there is a sequence of residues $i_1,\cdots,i_p$
and a sequence of multipartitions
$\uempty=\blam[0],\cdots,\blam[p]=\blam$ such that for
each $t$, $\blam[t-1]$ is obtained from $\blam[t]$ by removing all
the nodes in the ladder $i_t$-sequence of $\blam[t]$.
\end{dfn}

It is clear that a strong ladder multipartition is necessarily a
ladder multipartition. But it is not obvious that the converse is
also true. We conjecture that the converse is also true. Indeed we
shall prove that the converse is true in the case where $e=0$. Note
that in the case $r=2$, a strong ladder multipartition is the same
as a bipartition which has an optimal sequence in the sense of
\cite{AJ}.
\smallskip

The next definition is a natural generalization of \cite[Definition 4.1]{AJ}
to the arbitrary multipartition case.

\begin{dfn} \label{ladorder} Let $\blam,\bmu\in\mathcal{P}_r(n)$. We write $\blam\prec\bmu$ if there exist integers $1\leq s\leq r$ and
$t\geq 0$ such that \begin{enumerate}
\item $\lam^{(j)}=\mu^{(j)}$ for any $s+1\leq j\leq r$ and $\lam_j^{(s)}=\mu_j^{(s)}$ for any $j>t$;
\item $\lam_t^{(s)}<\mu_t^{(s)}$.
\end{enumerate}
\end{dfn}

It is clear that if $\bmu\lhd\blam$ then $\blam\prec\bmu$.

\begin{lem} \label{ladder} Any strong ladder multipartition is $(\bQ,e)$-restricted, and
hence is a Kleshchev multipartition with respect to $(e;\bQ)$.
\end{lem}

\begin{proof} Let $\blam=(\lam^{(1)},\cdots,\lam^{(r)})$ be a strong ladder multipartition in $\mathcal{P}_r(n)$.
By Definition \ref{ladm}, we can find a sequence of residues
$i_1,\cdots,i_p$ and a sequence of multipartitions
$\uempty=\blam[0],\cdots,\blam[p]=\blam$ such that for
each $t$, $\blam[t-1]$ is obtained from $\blam[t]$ by removing all
the nodes in the ladder $i_t$-sequence of $\blam[t]$. By the
definition of ladder sequence, we know that $i_s\neq i_t$ whenever
$|s-t|=1$. For each $1\leq t\leq p$, we use $a_t$ to denote the
number of semi-ladder $i_t$-nodes in the ladder $i_t$-sequence of
$\blam[t]$ and set $l_t:=\sum_{j=1}^{t}a_j$, $l_0=0$. Let
$\gamma_{l_{t-1}+1},\gamma_{l_{t-1}+2},\cdots,\gamma_{l_t}$ be the
ladder $i_t$-sequence of $\blam[t]$.

We claim that \begin{equation}\label{straighten}
f_{i_p}^{(a_p)}\cdots f_{i_1}^{(a_1)}\uempty=\blam+\sum_{\bmu\prec\blam}C_{\blam,\bmu}(v)\bmu,
\end{equation}
for some $C_{\blam,\bmu}(v)\in\mathbb{Z}_{\geq 0}[v,v^{-1}]$, where $f_{i_j}^{(a_j)}$ denotes the quantum dividing power
(cf. \cite[(1.4.1), (3.1.1)]{Lu}).

The proof is similar to the proof \cite[Proposition 4.2]{AJ} in the
bipartition case. We use induction on $p$. By definition, $l_p=n$.
Let $c:=n-a_p$. Then $\gamma_{c+1},\cdots,\gamma_{c+a_p}=\gamma_n$
is the ladder $i_p$-sequence for $\blam$. Let
$\blam':=\blam-\{\gamma_{c+1},\cdots,\gamma_{n}\}$, which is again a
strong ladder multipartition and hence $(\bQ,e)$-restricted, and
hence a Kleshchev multipartition with respect to $(e;\bQ)$.

By induction hypothesis, we have that $$
f_{i_{p-1}}^{(a_{p-1})}\cdots f_{i_1}^{(a_1)}\uempty=\blam'+\sum_{\bmu'\prec\blam'}C_{\blam',\bmu'}(v)\bmu'.
$$
Let $\bmu\neq\blam$ be a multipartition which appears in
$f_{i_p}^{(a_p)}\cdots f_{i_1}^{(a_1)} \uempty$ with
nonzero coefficient. Then there exist removable $i_p$-nodes
$\alpha_{1},\cdots, \alpha_{a_p}$ of $\bmu$ and a multipartition
$\bmu'$ such that $\bmu=\bmu'\sqcup\{\alpha_{1},\cdots,
\alpha_{a_p}\}$.

By the definition of ladder sequence of $\blam$, it is clear that $\bmu'=\blam'$ implies that $\bmu\prec\blam$. Hence we
can assume that $\bmu'\prec\blam'$. Suppose that $\bmu\nprec\blam$, i.e., $\blam\prec\bmu$. Then we can find integers $1\leq s\leq r$ and $t\geq 0$ such that \begin{enumerate}
\item ${\mu'}^{(l)}={\lam'}^{(l)}$ for any $s+1\leq l\leq r$ and ${\mu'}_j^{(s)}={\lam'}_j^{(s)}$ for any $j>t$;
\item ${\mu'}_t^{(s)}<{\lam'}_t^{(s)}$.
\end{enumerate}
We claim that
\begin{enumerate}
\item[c)] $\mu_j^{(l)}=\lam_j^{(l)}$ for any $(j,l)\in\{(j,l)|s+1\leq l\leq r\}\sqcup\{(j,s)|j>t\}$;
\item[d)] $\mu_{t+1}^{(s)}<\mu_{t}^{(s)}={\mu'}_{t}^{(s)}+1={\lam'}_{t}^{(s)}=\lam_t^{(s)}$;
\item[e)] at least one of the nodes $\gamma_{c+1},\cdots,\gamma_{c+a_p}$ is above $(t,\lam_t^{(s)},s)$.
\end{enumerate}
In fact, all of these statements follow from the fact that
$\bmu'\prec\blam'$, $\blam\prec\bmu$ and there are no addable
$i_p$-nodes below $\alpha_{a_p}$. Now d) implies that
$(t,\lam_t^{(s)},s)$ is an $i_p$-node of $\lam^{(s)}$. Hence it is
not a removable node of $\blam$ otherwise it has to be removed to
obtain $\blam'$ (by the definition of ladder $i_p$-sequence of
$\blam$). It follows that $\lam_{t+1}^{(s)}=\lam_t^{(s)}$. Thus
$\mu_{t+1}^{(s)}<\lam_t^{(s)}=\lam_{t+1}^{(s)}$, which is a
contradiction. This completes the proof of (\ref{straighten}).

Since $\bmu\prec\blam$ (i.e., $\blam\nprec\bmu$) implies that $\bmu\ntriangleleft\blam$, it follows from Lemma \ref{2dfn} that $\blam$ is $(\bQ,e)$-restricted. Applying Lemma \ref{mcor1}, we deduce that $\blam$ is a Kleshchev multipartition with respect to $(e;\bQ)$.
\end{proof}

Let $m$ be an arbitrary integer. We use $W_m$ to denote the
symmetric group generated by $s_i, i\in\Z/e\Z-\{m+e\Z\}$. Let
$W/W_m$ be the set of distinguished coset representatives of $W_m$
in $W$.

\begin{lem} {\rm (\!\!\cite[Lemma 3.2]{De})} \label{Dlem} Let $i\in\Z/e\Z$. Let $x\in W_m$, $d\in W/W_m$ and $w:=dx$. Suppose that
$s_iw<w$ and $s_id\not\in W/W_m$. Then $d^{-1}s_id=s_l$ for some $l\in\Z/e\Z$ with
$l\neq m+e\Z$ and such that $s_lx<x$.
\end{lem}

Recall that a multipartition $\blam=(\lam^{(1)},\cdots,\lam^{(r)})$
is said to be a multi-core if $\lam^{(j)}$ is an $e$-core for each
integer $1\leq j\leq r$. For each multipartition
$\blam=(\lam^{(1)},\cdots,\lam^{(r)})\in\mathcal{P}_r(n)$, we write
$$
\blam^{\diamond}:=(\lam^{(r)},\cdots,\lam^{(1)}).
$$
Clearly, $\blam^{\diamond}$ is a multi-core if and only if $\blam$
is a multi-core. The next theorem is the main result of this paper.

\begin{prop} \label{corecase} Let $\blam=(\lam^{(1)},\cdots,\lam^{(r)})\in\mathcal{P}_r(n)$. Suppose 
$\ulam^{\diamond}$ is a Kleshchev $r$-partition of $n$ with
respect to $(e,{v_r},\cdots,{v_1})$ and also a multi-core. Then
$\blam^{\diamond}$ is a strong ladder multipartition and hence
$(\bQ,e)$-restricted. Furthermore, for each ladder sequence
$\gamma_1>\cdots>\gamma_a$ of $\blam^{\diamond}$,
$\ulam^{\diamond}-\{\gamma_1,\cdots,\gamma_a\}$ is again a strong
ladder multipartition and hence a Kleshchev multipartition with
respect to $(e;{v_r},\cdots,{v_1})$.
\end{prop}
\begin{proof} We argue by induction on $n$. Suppose that the conclusion is true for any integer $0\leq n'<n$. In other words, for any multi-core Kleshchev $r$-partition $\blam'$ of $n'$
with respect to $(e;{v_r},\cdots,{v_1})$, $\blam'$ is a strong
ladder multipartition, and for any ladder sequence
$\gamma'_1>\cdots>\gamma'_{a'}$ of ${\blam'}$,
${\ulam'}-\{\gamma'_1,\cdots,\gamma'_{a'}\}$ is again a strong
ladder multipartition.

We now look at the multi-core Kleshchev $r$-partition
$\ulam^{\diamond}:=(\lam^{(r)},\cdots,\lam^{(1)})$ of $n$ with
respect to $(e;{v_r},\cdots,{v_1})$. By Lemma \ref{1ladlem}, there exists a ladder sequence $\gamma_1>\cdots>\gamma_a$
with $\res(\gamma_1)=i\in\Z/e\Z$. Suppose that the nodes
$\gamma_1>\cdots>\gamma_c$ are located in the component $\lam^{(t)}$
and $\gamma_{c+1}\not\in\lam^{(t)}$ for some $1\leq c<a$. Let $\umu$
be the $r$-partition which is obtained from $\ulam$ by deleting
the nodes $\{\gamma_1,\cdots,\gamma_c\}$. Since $\lam^{(t)}$ is an
$e$-core, $\lam^{(t)}$ has no addable $i$-nodes. By definition, $\gamma_1,\cdots,\gamma_c$ are all the removable 
$i$-nodes on the $e$-core partition $\lam^{(t)}$. It follows (by considering the abacus display of partition) that 
$\lam^{(t)}-\{\gamma_1,\cdots,\gamma_t\}$ is again an $e$-core, and hence 
$\umu^{\diamond}$ is again a multi-core. We are going to show that
$\umu^{\diamond}$ is Kleshchev. Note that
$\gamma_{c+1}>\cdots>\gamma_a$ is a ladder $i$-sequence of
$\umu^{\diamond}$. Once we can prove $\umu^{\diamond}$ is Kleshchev,
then by induction hypothesis that $$
\ulam^{\diamond}-\{\gamma_1,\cdots,\gamma_a\}=\umu^{\diamond}-\{\gamma_{c+1},\cdots,\gamma_{a}\}
$$ is a strong ladder multipartition, and then by definition, $\blam^{\diamond}$ must be a strong ladder
multipartition as well, which completes the proof of the proposition.

Let $\umu=(\mu^{(1)},\cdots,\mu^{(r)})$. By Theorem \ref{pthm} and
Lemma \ref{keylm}, we can find elements $w_1\geq w_2\geq\cdots\geq
w_r$ in $W$, such that $\lam^{(j)}=w_j\emptyset_{v_j}$, where the
subscript is used to indicate the charge ${v_j}$. There is a unique
way to write each $w_j$ in the form $d_jx_j$ where $d_j\in
W/W_{v_j}, x_j\in W_{v_j}$. For later use, we choose these elements
$w_j$ in a way such that $\sum_{j=1}^{r}\ell(x_j)$ is as small as
possible.\smallskip

Note that $\mu^{(j)}=\lam^{(j)}$ for any $j\neq t$. Since both $\lam^{(t)}$ and $\mu^{(t)}$ are $e$-cores, and $
\mu^{(t)}$ is obtained from $\lam^{(t)}$ by removing all its removable $i$-nodes, we deduce that $\mu^{(t)}=(s_id_t)\emptyset_{v_t}$ with
$d_t>s_id_t\in W/W_{v_t}$. Now for each integer $1\leq j\leq r$, we
define $$
w'_j=\begin{cases} w_j, &\text{if $j\neq t$;}\\
s_iw_j, &\text{if $j=t$.}
\end{cases}
$$
Note that $\ell(s_iw_t)=\ell(s_id_t)+\ell(x_t)=\ell(w_t)-1$. Hence
$w'_t<w_t$. Let $t+1\leq k\leq r$ be the smallest integer such that
$\lam^{(k)}$ contains addable $i$-nodes. Then for each integer
$t+1\leq l\leq k-1$, $\lam^{(l)}$ contains neither removable
$i$-nodes nor addable $i$-nodes.

We claim that \begin{equation} \label{equa30} w'_1\geq
w'_{2}\geq\cdots\geq w'_r.\end{equation} In fact, it suffices to
prove that $w'_t\geq w'_{t+1}$. Suppose this is not the case, then
we can deduce that $w_{t+1}$ must have a reduced expression which
starts from $s_i$ (otherwise the inequality $s_iw'_t=w_t\geq
w_{t+1}$ already implies that $w'_t\geq w_{t+1}$). Hence
$s_iw_{t+1}<w_{t+1}$. Note that since the $e$-core $\lam^{(t+1)}$ has no removable
$i$-nodes, $d_{t+1}$ has no reduced expression which starts from
$s_i$. It follows that $s_id_{t+1}>d_{t+1}$ and hence $s_id_{t+1}\not\in W/W_{v_{t+1}}$ 
(otherwise $s_iw_{t+1}=s_id_{t+1}x_{t+1}>d_{t+1}x_{t+1}=w_{t+1}$, a contradiction).
Applying Lemma \ref{Dlem}, we get that
$d_{t+1}^{-1}s_id_{t+1}=s_l$ for some $l\neq v_{t+1}+e\Z$ and such that
$\ell(s_lx_{t+1})<\ell(x_{t+1})$. In particular, we see that
\begin{equation}\label{equa31}
s_iw_{t+1}=d_{t+1}(s_lx_{t+1})<d_{t+1}x_{t+1}=w_{t+1}.
\end{equation}
For each integer $t+1\leq j\leq k-1$, we define $$
\widetilde{w}_j=\begin{cases} w_j, &\text{if $s_iw_j>w_j$;}\\
s_iw_j, &\text{if $s_iw_j<w_j$.}
\end{cases}
$$
We write $\widetilde{w}_j=\widetilde{d}_j\widetilde{x}_j$, where
$\widetilde{d}_j\in W/W_{v_j}, \widetilde{x}_j\in W_{v_j}$. Then
from (\ref{equa31}) we see that $\widetilde{w}_{t+1}=s_iw_{t+1}$,
$\widetilde{d}_{t+1}=d_{t+1}$ and
$\widetilde{x}_{t+1}=s_lx_{t+1}<x_{t+1}$. In general, for each
integer $t+1\leq j\leq k-1$, if $s_iw_j>w_j$, then by definition
$\widetilde{w}_j=w_j$, $\widetilde{d}_j=d_j$ and
$\widetilde{x}_j=x_j$; while if $s_iw_j<w_j$, then as $\lam^{(j)}$
has no removable $i$-nodes, it follows that
$s_id_j\not\in W/W_{v_j}$ (otherwise $s_iw_{j}=s_id_{j}x_{j}>d_{j}x_{j}=w_{j}$, a contradiction). Applying Lemma \ref{Dlem}, we get
that $d_{j}^{-1}s_id_{j}=s_l$ for some $l\neq v_{j}+e\Z$ and such that
$\widetilde{d}_{j}=d_{j}$, $\widetilde{x}_{j}=s_lx_j<x_{j}$. In
particular, we see that $\widetilde{w}_jW_{v_j}=w_jW_{v_j}$ for any
integer $t+1\leq j\leq k-1$. We claim that
\begin{equation}\label{equa32} w_1\geq\cdots\geq
w_t\geq\widetilde{w}_{t+1}\geq\cdots\geq\widetilde{w}_{k-1}\geq
w_{k}\geq\cdots\geq w_{r}.
\end{equation}
It is enough to show that
$\widetilde{w}_{t+1}\geq\cdots\geq\widetilde{w}_{k-1}\geq w_{k}$.
For each integer $t+1\leq j\leq k-2$, there are only the following
three possibilities: \medskip

\noindent {\it Case 1.} $\widetilde{w}_{j}=w_j$. In this case, it is
trivial to see that $\widetilde{w}_{j}\geq\widetilde{w}_{j+1}$.
\smallskip

\noindent {\it Case 2.} $\widetilde{w}_{j}=s_iw_j<w_j$,
$\widetilde{w}_{j+1}=s_iw_{j+1}<w_{j+1}$. From the inequality
$w_j\geq w_{j+1}$ it is also clear that
$\widetilde{w}_{j}\geq\widetilde{w}_{j+1}$.
\smallskip

\noindent {\it Case 3.} $\widetilde{w}_{j}=s_iw_j<w_j$,
$\widetilde{w}_{j+1}=w_{j+1}$. By definition, we know that
$s_iw_{j+1}>w_{j+1}$. In particular, $w_{j+1}$ has no
reduced expression starting from $s_i$. From the inequality
$s_i\widetilde{w}_j=w_{j}\geq w_{j+1}$ it follows that
$\widetilde{w}_{j}\geq w_{j+1}= \widetilde{w}_{j+1}$.
\smallskip

It remains to show that $\widetilde{w}_{k-1}\geq w_{k}$. If
$\widetilde{w}_{k-1}=w_{k-1}$, there is nothing to prove. Assume
$\widetilde{w}_{k-1}=s_{i}w_{k-1}<w_{k-1}$. Since
$\lam^{(k)}$ is an $e$-core which contains addable $i$-nodes. We
deduce that $d_{k}<s_id_{k}\in W/W_{v_{k}}$. In particular,
$\ell(s_iw_k)=\ell(s_id_kx_k)=\ell(s_id_k)+\ell(x_k)=\ell(w_k)+1$,
which implies that $w_k$ has no reduced expression which starts from
$s_i$. Therefore, from the inequality
$s_i\widetilde{w}_{k-1}=w_{k-1}\geq w_k$ we can deduce that
$\widetilde{w}_{k-1}\geq w_k$, as required. This completes the proof
of the claim (\ref{equa32}). \smallskip

Since $$
\sum_{j=t+1}^{k-1}\ell(\widetilde{x}_j)<\sum_{j=t+1}^{k-1}\ell({x}_j),
$$
we get a contradiction to our previous choice of $x_j$. Therefore,
we must have that $w'_t\geq w'_{t+1}$. This proves the claim (\ref{equa30}).
Now applying Theorem \ref{pthm}, we deduce
that $\umu^{\diamond}$ is Kleshchev, as required. This completes the proof of the proposition.
\end{proof}

\begin{cor} \label{mcor2} Let $\ulam:=(\lam^{(1)},\cdots,\lam^{(r)})$ be a Kleshchev $r$-partition of $n$ with respect to
$(e,{v_1},\cdots,{v_r})$. Let $\gamma$ be an arbitrary ladder node
of $\blam$, where $\res(\gamma)=i\in\Z/e\Z$. Suppose that $\ulam$ is
a multi-core. Then $\blam-\{\gamma\}$ is a again a Kleshchev
multipartition with respect to $(e;{v_1},\cdots,{v_r})$.
\end{cor}

\begin{proof} Let $\gamma=\gamma_1>\cdots>\gamma_a$ be the ladder $i$-sequence in $\blam$. It is clear that
$\bmu:=\blam-\{\gamma_1,\cdots,\gamma_{a}\}$ is again a multi-core.
By Proposition \ref{corecase}, we know that both $\blam$ and
$\bmu:=\blam-\{\gamma_1,\cdots,\gamma_{a}\}$ are strong ladder
multipartitions. Since $\{\gamma_2,\cdots,\gamma_a\}$ is the ladder
$i$-sequence of
$\blam-\{\gamma\}=\bmu\sqcup\{\gamma_2,\cdots,\gamma_a\}$, it
follows directly from definition that $\blam-\{\gamma\}$ is a strong
ladder multipartition. Now using Lemma \ref{ladder}, we see that
$\blam-\{\gamma\}$ must be a Kleshchev multipartition with respect
to $(e;{v_1},\cdots,{v_r})$ as well.
\end{proof}

To sum up, we have the following inclusion relations: $$\begin{aligned}
&\left\{\begin{matrix}\text{Strong ladder}\\
\text{$r$-partitions of $n$}\end{matrix}\right\}\subseteq\left\{\begin{matrix}\text{$(\bQ,e)$-restricted}\\
\text{$r$-partitions of $n$}\end{matrix}\right\}\subseteq\mathcal{K}_r(n),\\
&\left\{\begin{matrix}\text{Strong ladder}\\
\text{$r$-partitions of $n$}\end{matrix}\right\}\subseteq\left\{\begin{matrix}\text{Ladder}\\
\text{$r$-partitions of $n$}\end{matrix}\right\}.
\end{aligned}
$$
We conjecture they are actually all equalities. Proposition
\ref{corecase} says that
$$
\mathcal{K}_r(n)\bigcap\Bigl\{\text{multi-cores}\Bigr\}\subseteq\left\{\begin{matrix}\text{Strong ladder }\\
\text{$r$-partitions of $n$}\end{matrix}\right\}.
$$

In the remaining part of this section, we shall show that our
conjecture is indeed true in the case $e=0$. In particular, we shall
show that the ``only if" part of Conjecture \ref{GDJM} is always
true if $e=0$ and the notion of ladder multipartition coincides with
the notion of strong ladder multipartition in that case.

\begin{prop} \label{core2case} Suppose that $e=0$. Then any Kleshchev multipartition in $\mathcal{K}_r(n)$ is a strong
ladder multipartition and hence is $(\bQ,e)$-restricted. In that
case, for any ladder node $\gamma$ of a strong ladder multipartition
$\blam$, $\blam-\{\gamma\}$ is again a strong ladder multipartition.
\end{prop}
\begin{proof} Since in the case $e=0$, every multipartition is an
$e$-core. The proposition follows immediately from Proposition
\ref{corecase} and Corollary \ref{mcor2}.
\end{proof}

\begin{thm}  Suppose that $e=0$. Let $\blam\in\mathcal{P}_r(n)$. Then $\blam$ is a ladder multipartition if and only if $\blam$ is a strong ladder multipartition.
\end{thm}

\begin{proof} It suffices to show that if $\blam$ is a ladder multipartition, then $\blam$ is a strong ladder multipartition.

We make induction on $n$. By definition, $\blam$ has a ladder node
$\gamma$ such that $\blam-\{\gamma\}$ is again a ladder
multipartition. Write $\res(\gamma)=i\in\Z/e\Z$. Suppose that
$\gamma\in\lam^{(c)}$. Since $e=0$, $\gamma$ must be the unique
$i$-node of $[\lam^{(c)}]$. By induction hypothesis,
$\bmu:=\blam-\{\gamma\}$ is a strong ladder multipartition. In
particular, $\bmu\in\mathcal{K}_r(n-1)$. If $\bmu$ has no ladder
$i$-node, then $\gamma$ is already a ladder $i$-sequence of $\blam$.
In that case it follows from definition that $\blam$ is a strong
ladder multipartition. Therefore, we can assume that $\bmu$ does
have ladder $i$-nodes. Let $\gamma_1>\cdots>\gamma_a$ be the ladder
$i$-sequence of $\bmu$. By Proposition \ref{core2case},
$\bmu-\{\gamma_1,\cdots,\gamma_a\}$ is again a strong ladder
multipartition. Since $\gamma>\gamma_1>\cdots>\gamma_a$ is the
ladder $i$-sequence of
$\blam=\bmu\sqcup\{\gamma,\gamma_1,\cdots,\gamma_a\}$, it follows
directly from definition that $\blam$ must be a strong ladder
multipartition as well. This completes the proof of the theorem.
\end{proof}

\bigskip
\section*{Acknowledgments}

\thanks{The research was supported by an Australian Research
Council discovery grant and partly by the National Natural Science
Foundation of China (Project 10771014 \& 10871023). The author was
grateful to Professor Andrew Mathas for many helpful
discussion. He also thanks the referees for their helpful comments.}\bigskip


\begin{thebibliography}{10}

\bibitem{A1} {\sc S. Ariki}, {\em On the semi-simplicity of the
Hecke algebra of $(\mathbb{Z}/r\mathbb{Z})\wr {\mathfrak{S}_n}$}, J.
Alg., {\bf 169} (1994), 216--225.

\bibitem{A2} {\sc S. Ariki}, {\em On the decomposition
  numbers of the {Hecke} algebra of {$G(m,1,n)$}}, J.~Math. Kyoto Univ., {\bf
  36} (1996), 789--808.

\bibitem{A3}
{\sc S.~Ariki}, {\em On the
classification of simple modules for cyclotomic Hecke algebras of
type $G(m,1,n)$ and Kleshchev multi-partitions},  Osaka J. Math.,
(4) {\bf 38} (2001), 827--837.

\bibitem{A4} {\sc S. Ariki}, {\em Representations of
quantum algebras and combinatorics of Young tableaux}, Translated
from the 2000 Japanese edition and revised by the author, University
Lecture Series, {\bf 26}, American Mathematical Society, Providence,
RI, 2002.

\bibitem{AK}
{\sc S.~Ariki and K.~Koike}, {\em A {H}ecke algebra of
$\Z/r\Z\wr{\mathfrak{S}}\sb n$ and construction of its irreducible
  representations}, Adv. Math., {\bf 106} (1994), 216--243.

\bibitem{AJ}
{\sc S.~Ariki and N.~Jacon}, {\em Dipper--James--Murphy's conjecture
for Hecke algebras of type $B$}, preprint, arXiv: math/0703447, (2007).

\bibitem{AM} S.~Ariki and A.~Mathas, {\em The number of simple
modules of the Hecke algebras of type $G(r,1,n)$}, Math. Z., {\bf
233} (2000), 601--623.

\bibitem{AKT}
{\sc S.~Ariki, V.~Kreiman and S.~Tsuchioka}, {\em On the tensor
product of two basic representations of $U_v(\hat{sl}_e)$}, Adv.
Math., {\bf 218} (2008), 28--86.

\bibitem{BK} J.~Brundan and A.~Kleshchev, {\em Graded
  decomposition numbers for cyclotomic Hecke algebras}, Adv.
Math., {\bf 222} (2009), 1883--1942.

\bibitem{BM}
{\sc M.~Brou\'e and G.~Malle}, {\em {Zyklotomische Heckealgebren}},
Ast\'erisque, {\bf 212} (1993), 119--189.

\bibitem{De}
{\sc V.V.~Deodhar}, {\em Some characterizations of Bruhat ordering on a Coxeter group
and determination of the relative M\"obius function}, Invent. Math., {\bf 39} (1977),
187--198.

\bibitem{DJM}
{\sc R.~Dipper, G.~James, and A.~Mathas}, {\em Cyclotomic
$q$--{Schur} algebras}, Math.~Z., {\bf 229} (1999), 385--416.

\bibitem{DJM2}
{\sc R.~Dipper, G.~James, and E.~Murphy}, {\em Hecke algebras of
type $B_n$ at roots of unity}, Proc. London Math. Soc., (3) {\bf 70}
(1995), 505--528.

\bibitem{GL}
{\sc J.~J. Graham and G.~I. Lehrer}, {\em Cellular algebras},
Invent. Math., {\bf 123} (1996), 1--34.

\bibitem{Ka0} {\sc M.~Kashiwara}, {\em On crystal bases of the $q$-analogue of universal enveloping algebras},
Duke Math. J., {\bf 63}, (1991), 465--516.

\bibitem{Ka} {\sc M.~Kashiwara}, {\em Base cristallines des groups quantiques},
Cours Sp\'{e}c., vol. {\bf 9}, Soc. Math. France, 2002.


\bibitem{L1}
{\sc P.~Littelmann}, {\em A Littlewood--Richardson rule for symmetrizable Kac--Moody algebras},
Invent. Math., {\bf 116} (1994), 329--346.

\bibitem{L2}
{\sc P.~Littelmann}, {\em A plactic algebra for semisimple Lie
algebras}, Adv. Math., {\bf 124} (1996), 312--331.

\bibitem{Lu} {\sc G.~Lusztig}, {\em Introduction to Quantum Groups},
Progress in Math. {\bf 110} Birkh\"auser, Boston, 1990.

\bibitem{Ma1}
{\sc A.~Mathas}, {\em {Hecke algebras and Schur algebras of the
symmetric
  group}}, Univ. Lecture Notes, {\bf 15}, Amer. Math. Soc., 1999.

\bibitem{Ma2}
{\sc A.~Mathas}, {\em {Simple modules of Ariki-Koike algebras}}, Proc. Pure Symp. Math., {\bf 63}, (1998),
383--396.

\bibitem{MM}
{\sc T.~Misra and K.C.~Miwa}, {\em Crystal bases for the basic representations of $U_q(\widehat{\mathfrak{sl}}(n))$},
Comm. Math. Phys., {\bf 134} (1990), 79--88.



\end{thebibliography}

\bigskip\bigskip
\end{document}